\documentclass[12pt]{article}
\usepackage{amsmath}
\usepackage{amssymb} 
\begin{document}

\newtheorem{df}{Definition}
\newtheorem{thm}{Theorem}
\newtheorem{lm}{Lemma}
\newtheorem{pr}{Proposition}
\newtheorem{co}{Corollary}
\newtheorem{re}{Remark}
\newtheorem{note}{Note}
\newtheorem{claim}{Claim}

\begin{center}
\begin{huge}$\ell_p$ ($p>2$) Does Not Coarsely Embed Into a Hilbert Space\end{huge}

\indent

William~B.~JOHNSON and N.~Lovasoa~RANDRIANARIVONY \footnote{Both authors were supported in part by 
NSF 0200690 and Texas Advanced Research Program 010366-0033-20013.

This paper represents a portion of the second authors's dissertation being prepared at 
Texas A\&M University under the direction of the first author.}

\end{center}

\indent

A (not necessarily continuous) map $f$ between two metric spaces $(X,d)$ and $(Y, \delta )$ is called a 
{\sl coarse embedding\/} (see \cite[7.G]{G}) if there exist two non-decreasing functions $\varphi _1:[0,\infty)
\rightarrow [0,\infty )$ and  $\varphi _2:[0,\infty) \rightarrow [0,\infty )$ such that 
 
\begin{enumerate}
\item $\varphi_1(d(x,y))\leq \delta (f(x),f(y)) \leq \varphi_2(d(x,y))$
\item $\varphi _1(t) \rightarrow \infty$ as $t\rightarrow \infty$.
\end{enumerate}

Nowak \cite{N}, improving a theorem due to A. N. Dranishnikov, G. Gong, V. Lafforgue, and G. Yu \cite{DGLV}, 
gave a characterization of 
coarse embeddability of general  metric spaces into a Hilbert space using a result of Schoenberg on negative 
definite kernels.  He used
this characterization to show that the spaces $L_p(\mu)$ coarsely embed into a Hilbert space for $p<2$.  
In this article, we show that
$\ell_p$ does not coarsely embed into a Hilbert space when $p>2$.  It was already proved 
 in \cite{DGLV} that the Lipschitz universal space $c_0$ (see \cite{A}) does not
coarsely embed into a Hilbert space.
 
In its full generality, the statement of our result is as follows:

\indent

\begin{thm} Suppose that a Banach space $X$ has a normalized symmetric basis $(e_n)_n$ and that $\displaystyle
\liminf_{n
\rightarrow \infty} n^{-\frac{1}{2}}\|\sum_{i=1}^n e_i\|=0$.  Then $X$ does not coarsely embed into a Hilbert
space.\end{thm}

\indent

In \cite{Y} Yu proved that a discrete metric space with bounded geometry must satisfy the 
coarse geometric Novikov conjecture if it coarsely embeds into a Hilbert space, and in \cite{KY}
G.~Kasparov and Yu proved that to get the same conclusion it is sufficient that the metric space
coarsely embeds into a uniformly convex Banach space.  Our theorem suggests that the result of
\cite{KY} cannot be deduced from the earlier theorem in \cite{Y}, but as yet there is no example
of a discrete metric space with bounded geometry which coarsely embeds into $\ell_p$ for some
$2<p<\infty$ but not into $\ell_2$.  (The reader should be warned that what we called a ``coarse
embedding" is called a ``uniform embedding" is many places, including \cite{DGLV}, \cite{KY},
and \cite{Y}.  Following \cite{N}, we use the term coarse embedding to avoid confusion with the
closely related notion of uniform embedding as it is used in non linear Banach space theory
\cite{BL}; i.e., a bi-uniformly continuous mapping.)

Besides Schoenberg's classical work \cite{S} on positive definite functions, an important tool for 
proving the theorem is Theorem 5.2 in \cite{AMM}, which asserts that the hypothesis on $X$ in the 
theorem implies that every
symmetric continuous positive definite function on $X$ is constant.  We present the proof in five steps.

\indent

\begin{center}STEP 0: REDUCING TO THE $\alpha$-H\"OLDER CASE\end{center}

Let $f:X \rightarrow H$ be a coarse embedding satisfying

\begin{enumerate}
\item $\varphi_1(\|x-y\|)\leq \|f(x)-f(y)\| \leq \varphi_2(\|x-y\|)$
\item $\varphi _1(t) \rightarrow \infty$ as $t\rightarrow \infty$.
\end{enumerate}

Our first claim is that we do not lose generality by assuming that $\varphi_2(t)=t^{\alpha}$ with 
$0<\alpha<\frac{1}{2}$.

\indent

To prove this claim, note first that $(x,y) \mapsto \|f(x)-f(y)\|^2$ is a negative definite kernel on $X$.  
This can be seen by direct computations (see \cite[Proposition 3.1]{N}).  We refer the reader to 
\cite[Chapter 8]{BL} or
\cite[Section 2]{N} for the definitions of negative definite kernels and negative definite functions.  

So, (\cite[Lemma 4.2]{N}) for any $0<\alpha <1$, the kernel $N(x,y)= \|f(x)-f(y)\|^{2\alpha}$ is also negative definite and 
satisfies $N(x,x)=0$ (such a negative definite kernel is called normalized).

As a result, a theorem of Schoenberg (\cite{S} and   \cite[Chapter 8]{BL}) allows us to find a Hilbert space 
$H_{\alpha}$ and a function $f_{\alpha}:X \rightarrow H_{\alpha}$ such that $N(x,y)=\|f_{\alpha}(x)-f_{\alpha}(y)\|^2$.

\indent

On the other hand, since $X$, being a normed space, is (metrically) convex, the original function $f:X \rightarrow H$
is Lipschitz for large distances  (see e.g.   \cite[proof of Proposition 1.11]{BL}).  Consequently, without loss of
generality, we can assume by rescaling that we have the following for $\|x-y\|\geq 1$:

$$\|f(x)-f(y)\|_{_H} \leq \|x-y\|$$

\noindent and

$$\left ( \varphi_1(\|x-y\|) \right )^{\alpha} \leq \|f_{\alpha}(x)-f_{\alpha}(y)\|_{_{H_{\alpha}}} 
\leq \|x-y\|^{\alpha}$$

\indent

Now, let $\mathcal{N}$ be a 1-net in $X$ (i.e. $\mathcal{N}$ is a maximal 1-separated subset of $X$).  
The restriction of $f_{\alpha}$ to $\mathcal{N}$ is $\alpha$-H\"older, so if $0<\alpha<\frac{1}{2}$, then we can
extend $f_{\alpha}$ to an $\alpha$-H\"older map $\widetilde{f_{\alpha}}$ defined on the whole of $X$ (see \cite{WW},
last statement of Theorem 19.1):

$$\widetilde{f_{\alpha}}: X \rightarrow H_{\alpha}$$

$$\forall x \in \mathcal{N},~~ \widetilde{f_{\alpha}}(x)=f_{\alpha}(x)$$

\noindent and

$$\forall x,y \in X, ~~\|f_{\alpha}(x)-f_{\alpha}(y)\|_{_{H_{\alpha}}} \leq \|x-y\|^{\alpha}.$$

\indent

Now, write $\widetilde{\varphi}_1(t)=\inf \{\|\widetilde{f_{\alpha}}(x)-\widetilde{f_{\alpha}}(y)\|,~~\|x-y\|\geq t\}$,
 and observe that $\widetilde{\varphi}_1(t)\rightarrow \infty$ as $ t \rightarrow \infty$.  

\indent

This finishes the proof of our reduction to the case where $f$ is $\alpha$-H\"older and thus uniformly continuous.  
So from now on we will assume that our coarse embedding is a map $f:X \rightarrow H$ satisfying the following for all
$x,y \in X$:
$$\varphi_1(\|x-y\|)\leq \|f(x)-f(y)\| \leq \|x-y\|^{\alpha}$$

\noindent where $\varphi_1(t)\rightarrow \infty$ as $t \rightarrow \infty$.

\newpage

\begin{center}STEP 1\end{center}

Set $N(x,y)=\|f(x)-f(y)\|^2$.  Then $N$ is a normalized (i.e. $N(x,x)=0$) negative definite kernel on $X$, (see \cite[Proposition 3.1]{N}).  Now if we write $\phi_1(t)=(\varphi_1(t))^2$ and $\phi_2(t)=t^{2\alpha}$, then $N$ satisfies:

$
\begin{cases}
\phi_1(\|x-y\|)\leq N(x,y) \leq \phi_2(\|x-y\|), \\
\phi_1(t) \rightarrow \infty ~\text{as} ~t \rightarrow \infty .\\
\end{cases}
$

\indent

\begin{center}STEP 2\end{center}

The argument in this paragraph comes from \cite[Lemma 3.5.]{AMM}.

Let $\mu$ be an invariant mean on the bounded functions on $X$ (see e.g. \cite{BL} for the definition of invariant
means). 
 Define:

$$g(x)=\int_{X} N(y+x,y)\,d \mu(y)$$

Then we have the following for $g$:

\begin{itemize}
\item $g$ is well-defined because the map $y \mapsto N(y+x,y)$ is bounded for each $x\in X$

\item $\displaystyle g(0)=\int_{X} N(y,y)\,d \mu(y)=0$

\item For scalars $(c_i)_{1 \leq i \leq n}$ satisfying $\displaystyle \sum_{i=1}^{n} c_i =0$, we have:

\begin{equation*}
\begin{split}
\sum_{i,j=1}^n c_ic_jg(x_i-x_j)& =\displaystyle \sum_{i,j} c_ic_j\int_{X} N(y+x_i-x_j,y)\,d \mu(y) \\
 & \\
 & = \sum_{i,j=1}^n c_ic_j\int_{X} N(y+x_i,y+x_j)\,d \mu(y) \\
 & \\
 & = \int_{X} \left ( \sum_{i,j=1}^n c_ic_jN(y+x_i,y+x_j) \right )\,d \mu(y) \\
 & \\
 & = \int_{X} (\leq 0 )\,d \mu(y) \\
 & \\
 & \leq 0 \\
\end{split}
\end{equation*}

This is because $\mu$ is translation invariant, and $N$ is negative definite.  
This shows that $g$ is a negative definite function on $X$.

\item Finally, since $\displaystyle \int_{X} \,d \mu(y)=1$, we have:

$$\phi_1(\|x\|) \leq g(x) \leq \phi_2(\|x\|).$$

\end{itemize}

In summary, we have found  a negative definite function $g$ on $X$ which satisfies $g(0)=0$ and $\phi_1(\|x\|) \leq
g(x)
\leq
\phi_2(\|x\|)$, where $\phi_1(t) \rightarrow \infty$ as $t \rightarrow \infty$.

\indent

\begin{center}STEP 3\end{center}

Let $(e_n)_n$ be the normalized symmetric basis for $X$.  This means that for any choice 
of signs $(\theta_n)_n \in \{-1,+1\}$ and any choice of permutation $\sigma :\mathbb{N} \rightarrow \mathbb{N}$,

$$\|\sum_n \theta_n a_n e_{\sigma(n)}\|_{X}=\|\sum_n a_n e_n\|_{X}.$$

\indent

The purpose of this paragraph is to show that the negative definite function $g$ we 
found in the previous paragraph can be chosen to be symmetric, i.e. to satisfy for any choice of signs $(\theta_n)_n
\in \{-1,+1\}$ and any choice of permutation $\sigma :\mathbb{N} \rightarrow \mathbb{N}$ the equality:

$$g \left (\sum_n \theta_n a_n e_{\sigma(n)}\right )=g \left (\sum_n a_n e_n\right ).$$

For $x=\displaystyle \sum_{n=1}^{\infty} x_n e_n \in X$, define $g_m(x)$ to be the 
average of $\displaystyle g \left (\sum_{n=1}^{\infty}\theta_n x_n e_{\sigma (n)}\right )$ over all choices of signs
$\theta$ and permutations $\sigma$ with the restrictions that $\theta_n=1$ for $n>m$ and $\sigma(n)=n$ for $n>m$.

It follows that for all such $\theta, \sigma$, and for all $x=\displaystyle \sum_{n=1}^{\infty} x_n e_n \in X$,

$$\displaystyle g_m\left (\sum_{n=1}^{\infty} \theta_n x_n e_{\sigma (n)}\right )= 
\displaystyle g_m\left (\sum_{n=1}^{\infty} x_n e_n\right ).$$

Moreover, we also have 

$$\phi_1(\|x\|) \leq g_m(x) \leq \phi_2(\|x\|).$$

Next we show that the sequence $(g_m)_m$ is equicontinuous.  To check this, let us first check the 
continuity of $g$:

\begin{equation*}
\begin{split}
\|g(a)-g(b)\|
 & \leq \displaystyle \int_X \left |N(y+a,y)-N(y+b,y) \right | \,d \mu (y)\\
 & \\
 &= \displaystyle \int_X \left | \|f(y+a)-f(y)\|^2-\|f(y+b)-f(y)\|^2\right |\,d \mu (y)\\
 & \\
 &=\displaystyle \int_X \left ( \|f(y+a)-f(y)\|+\|f(y+b)-f(y)\|\right ) \cdot \\ & \phantom{XXX} 
\left | \|f(y+a)-f(y)\|-\|f(y+b)-f(y)\|\right |\,d \mu (y)\\
 & \\
 &\leq \displaystyle \int_X \left ( \|f(y+a)-f(y)\|+\|f(y+b)-f(y)\|\right ) \cdot \\
 & \phantom{XXX}\|f(y+a)-f(y+b)\|\,d \mu (y)\\
 & \\
 &\leq \displaystyle \int_X \left ( \|a\|^{\alpha}+\|b\|^{\alpha}\right )~\|a-b\|^{\alpha}\,d \mu (y).\\
\end{split}
\end{equation*}

So $\|g(a)-g(b)\| \leq \|a-b\|^{\alpha} \left ( \|a\|^{\alpha}+\|b\|^{\alpha}\right )$ and $g$ is continuous.

\indent

Now for the equicontinuity of $(g_m)_m$:

\indent

\begin{equation*}
\begin{split}
\|g_m(a)-g_m(b)\|
 &=\| \text{ave}~~\left ( g\left (\sum \theta_n a_n e_{\sigma(n)}\right ) - g\left 
(\sum \theta_n b_n e_{\sigma(n)}\right )\right )\| \\
 & \\
 & \leq \text{ave}~~ \|g\left (\sum \theta_n a_n e_{\sigma(n)}\right ) - g\left 
(\sum \theta_n b_n e_{\sigma(n)}\right )\| \\
 & \\
& \leq \text{ave}~~\left ( \|\sum \theta_n a_n e_{\sigma(n)}-\sum \theta_n b_n e_{\sigma(n)}\|^{\alpha} \right . \\
& \phantom{XXXXX}\left . \cdot \left (\|\sum \theta_n a_n e_{\sigma(n)}\|^{\alpha}+\|
\sum \theta_n b_n e_{\sigma(n)}\|^{\alpha}\right )\right )\\
 & \\
 &=\text{ave}~~\left (\|a-b\|^{\alpha} \left ( \|a\|^{\alpha}+\|b\|^{\alpha}\right ) \right )\\
& \\
 &=\|a-b\|^{\alpha} \left ( \|a\|^{\alpha}+\|b\|^{\alpha}\right ).\\
\end{split}
\end{equation*}

\indent

So by Ascoli's theorem, there is a subsequence $(g_{m_k})_k$ of $(g_m)_m$ which 
converges pointwise to a continuous function $\widetilde{g}$.  The property of the $g_m$'s implies that
$\widetilde{g}$ must necessarily be symmetric.  We have that $\widetilde{g}(0)=0$, and that $\phi_1(\|x\|) \leq
\widetilde{g}(x) \leq \phi_2(\|x\|)$.  Finally, as it is easily checked that the $g_m$'s are negative definite
functions, it also follows easily that $\widetilde{g}$ is a negative definite function.

\indent

\begin{center}STEP 4\end{center}

There is a relation between negative and positive definite kernels as given by a result of Schoenberg \cite{S};  see
also   \cite[Chapter 8]{BL}.  This result states that a kernel $N$ on $X$ is negative definite if and only if $e^{-tN}$
is positive definite for every $t>0$.

\indent

Since $\displaystyle \liminf_{n \rightarrow \infty} \frac{\|e_1+e_2+\cdots +e_n\|}{\sqrt{n}}=0$, 
and $\widetilde{f}=e^{-\widetilde{g}}$ is a symmetric continuous  positive definite function on $X$, we conclude by a
theorem of Aharoni, Maurey and Mityagin (see \cite[Theorem 5.2]{AMM}), that $\widetilde{f}$ is constant.

On the other hand, $\widetilde{f}(0)= e^{-\widetilde{g}(0)}=1$, while 
$0 \leq \widetilde{f}(x) \leq e^{-\phi_1(\|x\|)} \rightarrow 0$ as $\|x\| \rightarrow \infty$.  This gives a
contradiction and finishes the proof.

\bigskip

\bigskip

\newpage

\noindent    William B.~Johnson\newline
             Department Mathematics\newline
             Texas A\&M University\newline
             College Station, TX, USA\newline
             E-mail: johnson@math.tamu.edu

\bigskip  

\noindent    N.~Lovasoa Randrianarivony\newline
             Department Mathematics\newline
             Texas A\&M University\newline
             College Station, TX, USA\newline
             E-mail: nirina@math.tamu.edu


\begin{thebibliography}{9999}

\bibitem[A]{A} I.~Aharoni, 
Every separable metric space is Lipschitz equivalent to a subset of $c_0$,
Israel J. Math. 19 (1974), 284-291.

\bibitem[AMM]{AMM} I.~Aharoni, B.~Maurey, and  B.~S.~Mityagin, 
 Uniform embeddings of metric spaces and of
 Banach spaces into Hilbert spaces, 
Israel J. Math. 52(3) (1985), 251--265.



\bibitem[BL]{BL} Y.~Benyamini and  J.~Lindenstrauss, 
Geometric Nonlinear Functional Analysis,  
Coll. Pub. 48, Amer. Math. Soc., Providence, RI (2000).



\bibitem[DGLY]{DGLV} A.~N.~Dranishnikov, G.~Gong, V.~Lafforgue, and G.~Yu, 
Uniform embeddings 
into Hilbert spaces and a question of Gromov, 
Canad. Math. Bull.  45(1) (2002), 60--70. 



\bibitem[G]{G} M.~Gromov, Asymptotic invariants for infinite groups, Geometric Group 
Theory, 
G.~A.~Nibto and M.~A.~Roller, ed., Cambridge University Press (1993), 1--295. 



\bibitem[KY]{KY} G.~Kasparov and G.~Yu,
The coarse geometric Novikov conjecture and uniform convexity,
(preprint).


\bibitem[N]{N} P.~Nowak, 
Coarse embeddings of metric spaces into Banach  spaces, (preprint).



\bibitem[S]{S}  I.~J.~Schoenberg, 
Metric spaces and positive definite functions, 
Trans. Am. Math. Soc. 44 (1938), 522--536. 


 
\bibitem[WW]{WW} J.~H.~Wells, L.~R.~Williams, 
Embeddings and extensions in Analysis, 
Springer-Verlag (1975).

\bibitem[Y]{Y} G.~Yu,
The coarse Baum-Connes conjecture for spaces which admit a uniform embedding into Hilbert space,
Invent. Math. 139 (2000), no. 1, 201--240.


\end{thebibliography}
\end{document}